\documentclass[12pt]{article}

\usepackage[hidelinks]{hyperref}
\usepackage{epsfig}
\usepackage{latexsym}
\usepackage{caption}
\usepackage{amsmath}
\usepackage{amsfonts}
\usepackage{amssymb}
\usepackage{graphicx}
\usepackage{xcolor}
\usepackage{algorithmic}

\usepackage{color}

\makeatletter

\newcommand{\Rmnum}[1]{\expandafter\@slowromancap\romannumeral #1@}
\makeatother

\makeatletter
\let\@fnsymbol\@arabic
\makeatother

\begin{document}

\newtheorem{theorem}{Theorem}[section]
\newtheorem{observation}[theorem]{Observation}
\newtheorem{corollary}[theorem]{Corollary}
\newtheorem{algorithm}[theorem]{Algorithm}
\newtheorem{definition}[theorem]{Definition}
\newtheorem{guess}[theorem]{Conjecture}
\newtheorem{claim}[theorem]{Claim}
\newtheorem{problem}[theorem]{Problem}
\newtheorem{question}[theorem]{Question}
\newtheorem{lemma}[theorem]{Lemma}
\newtheorem{proposition}[theorem]{Proposition}
\newtheorem{fact}[theorem]{Fact}

\makeatletter
\newcommand\figcaption{\def\@captype{figure}\caption}
\newcommand\tabcaption{\def\@captype{table}\caption}
\makeatother

\newtheorem{acknowledgement}[theorem]{Acknowledgement}

\newtheorem{axiom}[theorem]{Axiom}
\newtheorem{case}[theorem]{Case}
\newtheorem{conclusion}[theorem]{Conclusion}
\newtheorem{condition}[theorem]{Condition}
\newtheorem{conjecture}[theorem]{Conjecture}
\newtheorem{criterion}[theorem]{Criterion}
\newtheorem{example}[theorem]{Example}
\newtheorem{exercise}[theorem]{Exercise}
\newtheorem{notation}[theorem]{Notation}
\newtheorem{solution}[theorem]{Solution}
\newtheorem{summary}[theorem]{Summary}

\newenvironment{proof}{\noindent {\bf
		Proof.}}{\rule{3mm}{3mm}\par\medskip}
\newcommand{\remark}{\medskip\par\noindent {\bf Remark.~~}}
\newcommand{\pp}{{\it p.}}
\newcommand{\de}{\em}

\newcommand{\qf}{Q({\cal F},s)}
\newcommand{\qff}{Q({\cal F}',s)}
\newcommand{\qfff}{Q({\cal F}'',s)}
\newcommand{\f}{{\cal F}}
\newcommand{\ff}{{\cal F}'}
\newcommand{\fff}{{\cal F}''}
\newcommand{\fs}{{\cal F},s}

\newcommand{\G}{\Gamma}
\newcommand{\g}{\gamma}
\newcommand{\wrt}{with respect to }

\newcommand{\ch}{\rm ch}
\newcommand{\ra}{{\rm ran}}
\newcommand{\co}{{\rm col}}
\newcommand{\sco}{{\rm scol}}
\newcommand{\wc}{{\rm wcol}}
\newcommand{\dc}{{\rm dcol}}
\newcommand{\ar}{{\rm arb}}
\newcommand{\ma}{{\rm mad}}
\newcommand{\di}{{\rm dist}}
\newcommand{\tw}{{\rm tw}}
\newcommand{\col}{{\rm col}}
\newcommand{\wcol}{{\rm wcol}}
\newcommand{\adm}{{\rm adm}}
\newcommand{\edp}[2]{#1^{[\natural #2]}}
\newcommand{\epp}[2]{#1^{\natural #2}}
\newcommand{\ind}{{\rm ind}}

\newcommand{\mad}{{\rm mad}}


\title{Generalized oriented chromatic numbers and oriented generalized coloring numbers of graphs }

\renewcommand{\thefootnote}{\arabic{footnote}}
\author{
	Wenwen Li \footnotemark\footnotemark[1]$~^,$\footnotemark[3]
	\and
	Jaroslav Ne{\v s}et{\v r}il \footnotemark[1]$~^,$\footnotemark[2]~~~~
	\and
	Daqing Yang \footnotemark[1]$~^,$\footnotemark[4]~~~
}

\footnotetext[1]{School of Mathematical Sciences,
	Zhejiang Normal University, Jinhua, Zhejiang 321004, China}
\footnotetext[2]{Computer Science Institute of Charles University (IUUK and ITI),
	Malostransk\'e n\'am. 25, 11800 Praha 1, Czech Republic. E-mail: {\small\texttt{nesetril@iuuk.mff.cuni.cz}}.} 	
\footnotetext[3]{E-mail: {\small1433828444\texttt{@qq.com}}.}
\footnotetext[4]{Grant number: NSFC, China 12271489.  E-mail: {\small\texttt{dyang@zjnu.edu.cn}}.}

\maketitle

\begin{abstract}


For a graph $G$ and a positive integer $q$, we bound the $q$-transitive oriented chromatic number $\chi_o^q(G)$ in terms of the oriented weak coloring number $\overrightarrow{\wcol}_{2q}(G)$.
For orientations with no directed cycle of length at most $2q$, we compare $q$-transitive oriented colorings with $\vec P_{2q}$-preserving oriented colorings.
We also show that the oriented strong $(2k-1)$-coloring number bounds the strong $k$-coloring number of $G$.
These  results characterize bounded expansion in terms of oriented generalized coloring  numbers, generalized oriented chromatic number and homomorphisms of acyclic orientations.

\end{abstract}


\section{Introduction}
Graph theory and combinatorics contain a wealth of problems and results concerning colorings and partitions. Motivated by earlier work on the oriented chromatic number \cite{HN-90,HNZ-96,HNZ-93,MSW,N-Z-96,S-s17}, we introduce and study various generalizations of this concept, which are called directed path preserving coloring,  $q$-transitive oriented chromatic number. The first requires a homomorphism to map the ends of each short
directed path to an arc. The second requires the homomorphism to be injective on each such path.
We compare the resulting chromatic numbers and bound them in terms of oriented generalized coloring numbers.
We then use these bounds to characterize graph classes of bounded expansion.

\subsection{Oriented coloring and homomorphism}

Let $G=(V,E)$ be a graph, an orientation of an edge $e=\{u,v\}\in E$ is one of the two possible ordered
pairs $uv$ or $vu$ called arcs.
An oriented graph $\vec{G}=(V,A)$ is obtained from an undirected graph $G=(V,E)$ by assigning to each edge $\{u,v\}\in E$ exactly one of the directed arcs $(u,v)$ or $(v,u)$. In this case, $\vec{G}$ is called an orientation of $G$, and $G$ is the underlying graph of $\vec{G}$.


Let  $\vec{G}$ be an
oriented graph, and let $c$ be a mapping from $V (\vec{G})$ to a
set of colors $X$. The mapping $c$ is an \emph{oriented coloring} of $\vec{G}$ if
\begin{itemize}
	\item[(1)] $c(u) \neq c(v)$ for every arc $uv$, and
	\item[(2)] $c(u) = c(x)$ implies $c(v) \neq  c(w)$ for every two arcs $uv$ and  $wx$ ($u$ and $x$ not necessarily distinct).
\end{itemize}
Observe that condition $(2)$ implies in particular that any two vertices linked
by a directed $2$-path, that is, a directed path of length $2$, must be assigned distinct colors.
The  {\em oriented chromatic number of an oriented
	graph $\vec{G}$}, denoted $\chi_o(\vec{G})$, is then defined as the least cardinality of a set $X$ such that $\vec{G}$ has an oriented coloring by the set $X$.
The {\em oriented chromatic number of
	an undirected graph $G$}, denoted $\chi_o(G)$, is then defined as the maximum of the oriented chromatic
numbers over all its orientations.

Chromatic and oriented chromatic numbers can be equivalently defined in terms of
graph homomorphisms. A \emph{homomorphism} $\varphi$ of a graph $G$ to a graph $H$, denoted $\varphi:
G \rightarrow H$, is an edge-preserving mapping from $V (G)$ to $V (H)$ ($\varphi(u)\varphi(v)$ is an edge in $H$ whenever $uv$ is an edge in $G$). Such a mapping $\varphi$ is called an \emph{$H$-coloring} of $G$. The
chromatic number of an undirected graph $G$ is then the least $k$ such that $G$ admits a
homomorphism to the complete graph $K_k$ on $k$ vertices. Similarly, the oriented chromatic
number of an oriented graph  $\vec{G}$ is the least $k$ such that there exists an oriented target
graph  $\vec{H}$ on $k$ vertices such that  $\vec{G}$ admits a  $\vec{H}$-coloring (homomorphisms of oriented graphs are defined as arc-preserving vertex mappings).

Homomorphisms of digraphs and
undirected graphs have been studied as a generalization of graph coloring (\cite{HN-90,HNZ-96,HNZ-93,MSW,N-Z-96,S-s17}). In particular, we know that classes of graphs with
bounded genus, bounded degree or bounded treewidth have bounded oriented
chromatic number \cite{KSZ-97,S-g97}. However, having bounded chromatic number is
not a sufficient condition for an undirected graph for having bounded oriented
chromatic number. For instance, the family of bipartite graphs has unbounded
oriented chromatic number \cite{KSZ-97,S-g97}.


\subsection{The q-transitive oriented coloring}


Let $q$ be a positive integer, the {\em $q$-transitive oriented coloring of an oriented graph $\vec{G}$} is a homomorphism $\varphi$ of $\vec{G}$ to an oriented graph $\vec{H}$, denoted $\varphi: V(\vec{G}) \rightarrow V(\vec{H})$, such that whenever $\vec{P}=u_{0}\dots u_{s}$ is a directed path with $1\leq s \le q$, then $\varphi(u_0)\varphi(u_s)$ is an arc in $\vec{H}$.
Such a mapping $\varphi$ is called a {\em $q$-transitive $\vec{H}$-coloring of $\vec{G}$.} For $q=1$, this is the ordinary coloring.

The {\em $q$-transitive oriented chromatic number of an oriented graph  $\vec{G}$,}   denoted $\chi^q_o(\vec{G})$, is the least $k$ such that there exists an  oriented target graph $\vec{H}$ on $k$ vertices such that  $\vec{G}$ admits a  $q$-transitive $\vec{H}$-coloring.
As before, the {\em  $q$-transitive oriented chromatic number of an undirected graph $G$}, denoted $\chi^q_o(G)$, is defined as the maximum of the $q$-transitive oriented chromatic numbers over its orientations $\vec{G}$, that admit an oriented $q$-transitive coloring.

Note that the (traditional) homomorphism $\varphi$ of an oriented graph $\vec{G}$ to an oriented graph $\vec{H}$ ($\varphi(u)\varphi(v)$ is an arc in $\vec{H}$ whenever $uv$ is an arc in $\vec{G}$) is a $q$-transitive oriented coloring with $q=1$; in this sense, the $q$-transitive oriented coloring and $q$-transitive
oriented chromatic number are generalizations of oriented coloring and oriented chromatic number.

The $q$-transitive oriented chromatic number  $\chi^q_o(\vec{G})$ 
can also be defined via the
existence of (generalized) oriented colorings: a $q$-transitive oriented  $k$-coloring of an oriented graph $\vec{G}$ is
a mapping $c: V (\vec{G})\rightarrow \{1, 2, \ldots , k\}$  such that:
\begin{enumerate}
	\item[(1)] if $\vec{P}=u_{0}\dots u_{s}$ is a directed path and $s \le q$, then $c(u_0) \neq c(u_s)$;
	\item[(2)] if  $\vec{P}=u_{0}\dots u_{s}$ and $\vec{Q}=v_{0}\dots v_{t}$ are two directed paths with the same direction and $s,t \le q$, then $c(v_t) = c(u_0)$ $\Rightarrow$ $c(u_s) \neq c(v_0)$ ($v_t$ and $u_0$ not necessarily
	distinct).
\end{enumerate}
The  $q$-transitive oriented chromatic number $\chi^q_o(\vec{G})$ is the minimum $k$ for which there exists an oriented $q$-transitive $k$-coloring of $\vec{G}$.

Note that, for any orientation $\vec{G}$ of an undirected graph $G$, if  $\chi^q_o(\vec{G})$ is well defined, then $\vec{G}$ has no  directed cycle of length at most $ 2q$.
For this reason,  in the definition of  $\chi^q_o(G)$, we only consider orientations  $\vec{G}$ of $G$ such that $\vec{G}$ does not have any directed cycle of length at most $2q$, this restriction is understood whenever we state that $\chi^q_o(G)$ is bounded.


For any oriented graph $\vec{G}$ without directed cycles of length at most $2q$, define the $q$-th power of $\vec{G}$ to be the digraph $\vec{G^{q}}=(V', A')$, where $V'=V(\vec{G})$ and $A'=\{ (x, y):$   there exists a directed path from $x$ to $y$ which has length at most $q$ in $\vec{G}\}$. By the definition of $\vec{G^{q}}$, $\chi_{o}(\vec{G}^{q})$ is well defined and $\chi_{o}(\vec{G}^{q})=\chi^q_o(\vec{G})$.
First, we prove that $\chi_{o}(\vec{G}^{q})\geq \chi^q_o(\vec{G})$.
Suppose that $c$ is an oriented coloring of $\vec{G}^{q}$. For any directed path  $\vec{P}=u_{0}\dots u_{s}$ and $s\leq q$,
then $(u_{0}, u_{s})\in A(\vec{G}^{q})$, which implies that $c(u_{0})\neq c(u_{s})$;
for any two directed paths of length at most $q$ in $\vec{G}$,  $\vec{P}=u_{0}\dots u_{s}$, $\vec{Q}=v_{0}\dots v_{t}$ and $s,t \le q$, then $(u_{0}, u_{s}), (v_{0}, v_{t})\in A(\vec{G}^{q})$, so $c(v_t) = c(u_0)$ $\Rightarrow$ $c(u_s) \neq c(v_0)$. Conversely, one can verify that  $\chi_{o}(\vec{G}^{q})\leq \chi^q_o(\vec{G})$ also holds.



This paper is devoted to the study of $q$-transitive oriented coloring of graphs.
The newly introduced parameter $\chi^q_o(G)$ is naturally related to the concept of $\vec{T}$-preserving homomorphisms of graphs, where $\vec{T}$ is an oriented graph,  which was first introduced and studied by  Ne{\v s}et{\v r}il,  Sopena, and Vignal in \cite{NSV-97}.
In Section $2$, we introduce the concept of $\vec{T}$-preserving homomorphisms of graphs, and
study its relationship with the $q$-transitive oriented chromatic number. 
\subsection{The directed path preserving coloring}

\label{sec-w2-o}

%

 Observe that if $\varphi$ is a homomorphism of an oriented graph $\vec{G}$ to an oriented graph $\vec{H}$  then for every directed $2$-path $uvw$ in $\vec{G}$, $\varphi(u)\varphi(v)\varphi(w)$ is a directed $2$-path in $\vec{H}$. We say that $\varphi$ is $\vec{P}_2$-preserving.

It is then natural to generalize this notion as follows: let  $\vec{T}$ be any oriented graph; a homomorphism of an oriented graph
$\vec{G}$ to an oriented graph $\vec{H}$ is said to be $\vec{T}$-preserving if every connected subgraph
$\vec{G'}$ of $\vec{G}$  which is isomorphic to a subgraph of $\vec{T}$ is isomorphic to its image under
$\varphi$ in $\vec{H}$  (in other words, any two vertices of $\vec{G'}$  must be assigned distinct colors).
A {\em $\vec{T}$-preserving oriented coloring} of an oriented graph $\vec{G}$ will thus be an oriented
coloring such that any two vertices belonging to some subgraph isomorphic to a
subgraph of $\vec{T}$ are assigned distinct colors. We define the $\vec{T}$-preserving oriented chromatic number $\vec{\chi}_{\vec{T}} (G)$ of an oriented graph $\vec{G}$ as the minimum number
of vertices in an oriented graph $\vec{H}$  such that there exists a $\vec{T}$-preserving homomorphism of $\vec{G}$ to $\vec{H}$. As before, the $\vec{T}$-preserving oriented chromatic number
of an undirected graph is defined as the maximum of the $\vec{T}$-preserving oriented
chromatic numbers of its orientations.

This notion of preservation can also be considered in the undirected case.
However, families of undirected graphs with bounded $T$-preserving chromatic number (for nontrivial $T$) must have bounded degree. To see this, let $P_2$
denote the undirected $2$-path and $S_n$ the star with $n$ branches (that is
$V (S_n) = \{x, y_1, y_2, \ldots, y_n\}$ and $E(S_n) = \{xy_i, 1 \le i \le n\})$. Observe that if $S_n$ admits a $P_2$-preserving homomorphism to an undirected graph $G$ then $G$
cannot have fewer than $n + 1$ vertices. In other words, the $P_2$-preserving chromatic number of stars is unbounded. By contrast, every oriented star has a
$\vec{P}_2$-preserving homomorphism to $\vec{P}_2$ itself and has $\vec{P}_2$-preserving oriented chromatic number $3$.

Note that $P_2$-preserving homomorphisms
of undirected graphs were called {\em local homomorphisms} \cite{N-71}.
Ne{\v s}et{\v r}il proved in \cite{N-71} that for every connected undirected graph $G$,
automorphisms are the only $P_2$-preserving homomorphisms of $G$ to itself. The same conclusion holds for $T$-preserving homomorphisms for every connected graph
$T$ with at least two edges.
Therefore, $\vec{T}$-preserving homomorphisms can be studied for a wider class of graphs in the oriented setting.
For directed graphs the situation is more complicated.
For instance, for every directed graph $\vec{T}$ there exists a
directed graph  $\vec{G}_{\vec{T}}$ with a non-injective $\vec{T}$-preserving homomorphism of  $\vec{G}_{\vec{T}}$  to
itself.

In \cite{NSV-97} it was observed that only families of graphs with bounded degree can have bounded $\vec{T}$-preserving oriented chromatic number when $\vec{T}$  has both in-degree and out-degree at least two. They also provide sufficient conditions under which families of oriented graphs have bounded  $\vec{T}$-preserving oriented chromatic number when $\vec{T}$ is a directed path or a directed tree.


\subsection{Bounded expansion and directed bounded expansion}

Bounded expansion and nowhere dense graph classes were
introduced by Ne\v{s}et\v{r}il and Ossona de Mendez  \cite{Book-N-O,N-O-I,N-O-II-1,N-O-II-2,N-O-II-3,N-O-II-4,N-O-III} to capture structural sparsity in undirected graphs; they encompass a wide range of uniformly sparse graph classes.
An important feature of these classes is that they admit many equivalent characterizations:
by the density of bounded-depth minors, by low tree-depth colourings \cite{N-O-II-1}, by generalised
colouring numbers \cite{Zhu-BGCN}, see also \cite{S.S2026} and many more.
Each of these different aspects of the theory comes with its own set of combinatorial tools.

Analogous to the corresponding notions for undirected graphs,  Kreutzer and Tazari  first introduced the notion of directed graph minors, the greatest reduced average density (grad) and  bounded expansion for digraphs in \cite{KT}.

 A digraph $\vec{H}$ has a $directed$ $model$ in a digraph $\vec{G}$ if there is a function $\delta$ mapping vertices $v \in V (\vec{H})$ of $\vec{H}$   to subgraphs $\delta(v)\subseteq \vec{G}$ and arcs $a \in A(\vec{H})$ to arcs $\delta(a) \in A(\vec{G})$ such that

1. if $v \neq u$, then $\delta(v) \cap \delta(u) = \emptyset$;

2. if $a = (u, v)$ and $\delta(a) = (u_{0}, v_{0})$ then $u_{0} \in \delta(u)$ and $v_{0}\in \delta(v)$. For $v \in V (\vec{H})$, let $in(\delta(v)):=V (\delta(v)) \cap \bigcup_{a=(u,v)\in A(\vec{H})} V(\delta(a))$ and $out(\delta(v)):= V(\delta(v)) \cap \bigcup_{a=(v,w)\in A(\vec{H})}V(\delta(a))$;

3. we require that for every $v \in V(\vec{H})$
(i) there is a directed path in $\delta(v)$ from every $u \in in(\delta(v))$ to every $u'\in out(\delta(v))$;
(ii) there is at least one source vertex $s_{v} \in \delta(v)$ that reaches (by a directed path in $\delta(v)$) every element of $out(\delta(v))$;
(iii) there is at least one sink vertex $t_{v} \in \delta(v)$ that can be reached (by a directed path in $\delta(v)$) from every element of $in(\delta(v)$).

We write $\vec{H}\preccurlyeq \vec{G}$ if $\vec{H}$ has a directed model in $\vec{G}$ and call $\vec{H}$ a $directed$ $minor$ of $\vec{G}$. For $r\geq 0$, a digraph $\vec{H}$ is a $depth$-$r$ $minor$ of a digraph $\vec{G}$, denoted as $\vec{H}\preccurlyeq _{r}\vec{G}$, if there exists a directed model $\delta$ of $\vec{H}$ in $\vec{G}$ in which for all $v\in V(\vec{H})$ and all $s\in in(\delta(v))$ and $t\in out(\delta{(v)})$ there is a path from $s$ to $t$ in $\delta(v)$ of length $\leq r$.

A digraph $\vec{H}$ is a $topological$ $minor$ of a digraph $\vec{G}$ if there is an injective function $\delta$ mapping vertices $v \in V (\vec{H})$ to vertices of $V(\vec{G})$ and arcs $a\in A(\vec{H})$ to directed paths in $\vec{G}$ such that if $a = (u, v) \in A(\vec{H})$, then $\delta(a)$ is a path from $\delta(u)$ to $\delta(v)$ in $\vec{G}$ which is internally vertex disjoint from all vertices $\delta(w)$ (for $w \in V(\vec{H})$) and all paths $\delta(a')$ (for $a'\in A(\vec{H})$, $a' \neq a$). For $r \geq 0$, $\vec{H}$ is a $topological$ $depth$-$r$ $minor$ of $\vec{G}$, written $\vec{H} \preccurlyeq^{top}_{r}\vec{G}$, if it is a topological minor and all paths $\delta(a)$ have length at most $2r$.

Let $\vec{G}$ be a digraph and let $r \geq 0$. The $greatest$ $reduced$ $average$ $density$ $of$ $rank$ $r$ (short $grad$) of $\vec{G}$ is
$$\nabla_{r}(\vec{G}):=\max \left\{\frac{|E(\vec{H})|}{|V(\vec{H})|}:\vec{H}\preccurlyeq_{r}\vec{G}\right\},$$
and its $topological$ $greatest$ $average$ $density$ $of$ $rank$ $r$ (short $top$-$grad$) is
$$\widetilde{\nabla}_{r}(\vec{G}):=\max \left\{\frac{|E(\vec{H})|}{|V(\vec{H})|}:\vec{H}\preccurlyeq^{top}_{r}\vec{G}\right\}.$$

\begin{definition}
A class $\mathcal{D}$ of digraphs has bounded expansion if there is a function $f$: $\mathbb{N}\rightarrow\mathbb{N}$ such that for all $r\geq0$ we have $\nabla_{r}(\vec{G})\leq f(r)$ (or equivalently, $\widetilde{\nabla}_{r}(\vec{G})\leq f(r)$ ) for all $\vec{G}\in \mathcal{D}$.
\end{definition}

In Corollary~\ref{undirected-directed-bounded-expension}, we show that an undirected graph class has bounded expansion if and only if the class of all orientations of its graphs
has directed bounded expansion.

In Section $3$, we introduce the relevant notions concerning bounded expansion classes and the various coloring parameters used in this paper.
In particular, we introduce directed versions of generalized coloring numbers.

In Section 4, we use  oriented generalized coloring numbers to derive an upper bound for $q$-transitive oriented chromatic number (our main motivation). The proof of this theorem is in Section \ref{S-q-dis-oc}.
Our proofs are based on the analysis of weak $k$-coloring numbers (initiated in \cite{KY,Zhu-BGCN}). This approach is particularly suitable in the present context, as demonstrated by the improvements we obtain over some earlier results and bounds.

In Section 5, we study the relationship between the generalized coloring numbers  and oriented generalized coloring numbers.

\section{Directed path preserving coloring and q-transitive  oriented  coloring}
\label{sec-w2-o-sub}

Proposition \ref{pro-oc-Pq} establishes a relation between $\chi^q_o(G)$ and the $\vec{P}_q$-preserving oriented chromatic number of an orientation of $G$.
A relation in the reverse direction is presented in Proposition \ref{pro-Pq-to-Oc}.

\begin{proposition} \label{pro-oc-Pq}
If $\chi^q_o(G) \le s$, then for any orientation  $\vec{G}$ of $G$ and $\vec{G}$ does not have any directed cycle of length at most $ 2q$,
we have $\vec{\chi}_{\vec{P}_{2q}}(\vec{G}) \le s$.
\end{proposition}
\begin{proof}
For an orientation  $\vec{G}$ of $G$ that does not contain a directed cycle of length at most $2q$,
suppose the oriented $q$-transitive $k$-coloring $c$ of $\vec{G}$ is
a mapping $c$ that witnesses  $\chi^q_o(\vec{G}) \le s$. It suffices to show that for any directed path $\vec{P}=u_{0}\dots u_{s}$ of length $s \le 2q$, all the $u_i$, $i\leq s$, are colored differently under $c$. Suppose, for a contradiction, that $c(u_i) = c(u_j)$, and $i<j$.
Since the two ends of any directed path of length at most $q$ are colored differently under $c$ (by (1) of the definition of $c$), we can assume that $j-i > q$.
Then  $\vec{P}' = u_{i} \dots u_{i + \lceil \frac{j-i}{2}\rceil }$ and $\vec{P}'' = u_{i + \lceil \frac{j-i}{2}\rceil } \dots u_{j}  $ are directed paths of length  $1 \le \lceil \frac{j-i}{2}\rceil \le q$;  by (2) of the definition of $c$, $c(u_i) \neq c(u_j)$, a contradiction.
\end{proof}

\begin{proposition} \label{pro-Pq-to-Oc}
Let $G$ be a graph. For any orientation  $\vec{G}$ of $G$ such that $\vec{G}$ does not have any directed cycle of length at most $2q$,
and any integer $q \ge 1$, if  $\vec{\chi}_{\vec{P}_{2q}} (\vec{G}) \le s$,  then $\chi^q_o(\vec{G}) \le  c_q(s) \le 2^{s}-1$.
\end{proposition}

\begin{proof}
For an orientation  $\vec{G}$ of graph $G$ and $\vec{G}$ does not have any directed cycle of length at most $2q$, suppose that a coloring $\theta$ witnesses $\vec{\chi}_{\vec{P}_{2q}} (\vec{G}) \le s$.
For any $u$, $v\in V(\vec{G})$, then $\theta(u)\neq\theta(v)$ if there exists a directed path between $u$ and $v$ which has length at most $2q$ in $\vec{G}$.

Let $U_{\vec{G^{q}}}$ be the undirected graph defined by $V(U_{\vec{G^{q}}})=V(\vec{G})$ and $uv$ is an edge in $U_{\vec{G}^{q}}$ if and only if $u$ and $v$ are joined in $\vec{G}^{q}$
by a directed path of length at most $2$. By the above, it follows that $\theta$ is a proper coloring of $U_{{\vec{G}}^{q}}$.


Let $\vec{G}^{q}$ be the $q$-th power of $\vec{G}$ and $\vec{H}_{s}$ be the oriented graph defined as follows: the vertices of $\vec{H}_{s}$ are all the tuples of the form $(i; a_{1}, \cdots, a_{i-1})$ with $1\leq i\leq s$ ($i$ is called the identity of the vertex) and $a_{j}\in \{0, 1\}$ for every $j$, $1\leq j\leq i-1$.
Let $a=(i; a_{1},\cdots, a_{i-1})$ and $b=(j; b_{1},\cdots, b_{j-1})$ be two vertices in $\vec{H}_{s}$ with $i< j$.
There is  an arc from $a$ to $b$ if $b_{i}=0$ and an arc from $b$ to $a$ otherwise.
Observe that the graph $\vec{H}_{s}$ has $2^{s}-1$ vertices and all arcs in $\vec{H}_{s}$ link vertices with distinct identities.

We denote the set of vertices $x$ with $\theta(x)=i$ by $V_{i}$ and the (oriented) subgraph of $\vec{G}^{q}$ induced by $V_{i}\cup V_{j}$ by $\vec{G}^{q}_{i,j}$.
From the definition of $U_{{\vec{G}}^{q}}$ we get that ${\vec{G}}^{q}_{i, j}$ contains no directed 2-path in ${\vec{G}^{q}}$, i. e.,
if $x\in V_{i}$ and $x\rightarrow V_{j}$, then there is no vertex $y\in V_{j}$ such that $y\rightarrow x$.
We then define a mapping
$\phi: V({\vec{G}}^{q})$ to $V(\vec{H}_{s})$, and for every $x\in V_{i}, 1\leq i\leq s$, let $\phi(x)=(i; a_{1},\cdots, a_{i-1})$ where for every $j$, $1\leq j< i$, $a_{j}=0$ if there is an arc from $V_{j} \rightarrow$ $x$ in ${\vec{G}}^{q}$ and $a_{j}=1$ otherwise. This implies that $\phi$ is an oriented $\vec{H}_{s}$-coloring of ${\vec{G}}^{q}$.
By the preceding construction, $\phi$ is clearly a $q$-transitive $\vec{H}_{s}$-coloring of $\vec{G}$. This completes the proof.

\end{proof}

\section{Generalized coloring  and  oriented generalized coloring}

\subsection{Generalized coloring of graphs}

It is well known that the chromatic number can be usefully studied relative to a vertex ordering, which provides an online setting.

Fix a graph $G=(V,E)$ and let $\Pi:=\Pi(V)$ be the set of total orderings of the vertex set $V$. For $\sigma\in\Pi$ and $x\in V$, set
\begin{enumerate}
\item $V_{\sigma}^{l}(x)=\{y\in V: y<_{\sigma}x\}$, $V_{\sigma}^{l}[x]=V_{\sigma}^{l}(x) \cup \{x\}$;
and
\item $V_{\sigma}^{r}(x)=\{y\in V: x<_{\sigma}y \}$, $V_{\sigma}^{r}[x]=V_{\sigma}^{r}(x)\cup \{x\}$.
\end{enumerate}

Thus $\{V_{\sigma}^{l}(x),\{x\},V_{\sigma}^{r}(x)\}$ partitions $V$
into the \emph{left set} of $x$, singleton $x$, and the \emph{right
set} of $x$. The coloring number of $G$, denoted $\co(G)$, is defined
by
\[
\co(G)=\min_{\sigma\in\Pi}\max_{x\in V}|N[x]\cap V_{\sigma}^{l}[x]|.
\]
A graph $G$ is $k$-degenerate if every nonempty subgraph of $G$ has a vertex of degree at most $k$. The coloring number of $G$
is one more than its degeneracy; under an older definition of degeneracy they were the same.
Greedily coloring the vertices of $G$ in an order that witnesses its coloring number shows that $\chi(G)\leq\co(G)$ where $\chi(G)$ denotes the chromatic number of $G$.

Our bounds will be stated in terms of generalized \emph{coloring numbers}.
These numbers were first introduced in \cite{KY} after similar notions were explored by various authors \cite{MR1198403,45,K-Comp,KT-OGC} in the cases $k=2,4$.

Let $k$ be a positive integer,  a vertex $y$ is \emph{weakly $k$-reachable} from $x$ with respect to $\sigma$ if $y\in V_{\sigma}^{l}[x]$ and there is an $x,y$-path $P$ of length at most $k$ whose vertices are contained in $V_{\sigma}^{r}[y]$; it is \emph{strongly $k$-reachable} from $x$ with respect to $\sigma$
if in addition the vertices of $P$ are contained in $V_{\sigma}^{r}[x]+y$.
Let $W_{\sigma}^{k}(x)$ be the set of vertices that are weakly $k$-reachable from $x$ with respect to $\sigma$ and $S_{\sigma}^{k}(x)$ be the set of vertices that are strongly $k$-reachable from $x$ with respect to $\sigma$.
The \emph{weak $k$-coloring number} of $G$, denoted $\wc_{k}(G)$, and the \emph{strong $k$-coloring number}
of $G$, denoted by $\co_{k}(G)$, are defined by:
\[
\wc_{k}(G)=\min_{\sigma\in\Pi}\max_{x\in V}|W_{\sigma}^{k}(x)|\text{ and }\co_{k}(G)=\min_{\sigma\in\Pi}\max_{x\in V}|S_{\sigma}^{k}(x)|.
\]

The following lemma is given in \cite{KY, Book-N-O}. 

\begin{lemma} \label{kc-L1}
Every graph $G$ satisfies ${\wcol}_{k}\left(G\right) \leq \left( {\col}_{k}\left( G\right) \right) ^{k}$.
\end{lemma}

Generalized coloring numbers are particularly useful in the study of structured graph families.
Ne{\v s}et{\v r}il and Ossona de Mendez introduced the class of graphs of \emph{bounded expansion} \cite{N-O-I,N-O-II-1,N-O-II-2,N-O-II-3,N-O-III}; a thorough treatment appears in  \cite{Book-N-O}. Zhu \cite{Zhu-BGCN} and  Yang \cite{Yang-GTF} proved an equivalent characterization via strong coloring numbers, which we adopt as our working definition.

\begin{definition}
A graph class $\mathcal C$ has bounded expansion if there exists a function $f:\mathbb Z^+\to \mathbb Z^+$ such that for all $G\in\mathcal C$ and all $k\in\mathbb N$,
\[
\mathrm{col}_k(G)\le f(k).
\]
\end{definition}
It is well known that every minor-closed graph class has bounded expansion \cite{Book-N-O}.

\subsection{Oriented generalized coloring of graphs}

\label{S-directed-col}

In this section, we introduce the concepts of the oriented $k$-coloring number and the oriented weak $k$-coloring
number of an (undirected) graph $G$.  These concepts are motivated by their undirected counterparts in \cite{KY}, and by the preliminary versions for $k\le2$ in \cite{Yang-Zhu-GCD}.

Previous work in \cite{KRSW} establishes relations among the directed parameters $\overrightarrow{\mathrm{adm}}_r$, $\overrightarrow{\mathrm{col}}_r$ and $\overrightarrow{\mathrm{wcol}}_r$, thereby establishing an equivalence between directed bounded expansion and the uniform boundedness of $\overrightarrow{\mathrm{wcol}}_r(\vec{G})$ by a function of $r$. Moreover, Kreutzer et al. \cite{KMORS} provide a polynomial-time vertex-ordering algorithm for digraph classes with this property.

Let $\vec{G}=(V,A)$ be a digraph, and let $\Pi:=\Pi(V(\vec{G}))$ be the set of total
orderings of the vertex set $V(\vec{G})$.
For $\sigma\in\Pi$ and $x, y \in V$, we say that $x$ is \emph{weakly } $(k,+)$\emph{-reachable from} $y$
in $\vec{G}$ if $x\le_{\sigma}y$ and there exists a directed path $\vec{P}$ from $y$ to  $x$ of length at most $k$ (i.e., with at most $k$ arcs)
such that every vertex $z$ of $\vec{P}$ satisfies $x \le_{\sigma}z$.
If in addition every internal vertex $z$ of  $\vec{P}$ satisfies $y<_{\sigma}z$ then $x$ is $(k,+)$\emph{-(strongly) reachable from} $y$.
Let $S_{\vec{G}, \sigma}^{k,+}(y)$ be the set of vertices that are $(k,+)$-reachable from $y$,
and let $W_{\vec{G}, \sigma}^{k,+}(y)$  be the set of vertices that are weakly $(k,+)$-reachable from $y$.
Note that by this definition, $y \in S_{\vec{G}, \sigma}^{k,+}(y)$ and $y \in W_{\vec{G}, \sigma}^{k,+}(y)$.

Define the $(k,+)$-\emph{coloring number} of $\vec{G}$, denoted by $\overrightarrow{\col}^+_{k}\left( \vec{G} \right) $, and the \emph{weak }$(k,+)$\emph{-coloring number} of $\vec{G}$, denoted by $\overrightarrow{\wcol}^+_{k}\left( \vec{G}\right) $, by
\begin{eqnarray*}
\overrightarrow{\col}^+_{k}\left( \vec{G}\right) &=&\min_{\sigma\in \Pi }
\max_{y\in V\left( \vec{G} \right) }\left|  S_{\vec{G}, \sigma}^{k,+}(y) \right|
\text{ and} \\
\overrightarrow{\wcol}^+_{k}\left(\vec{G}\right) &=&\min_{\sigma\in \Pi}
\max_{y\in V\left( \vec{G} \right) }\left|  W_{\vec{G}, \sigma}^{k,+}(y) \right|;
\end{eqnarray*}

Analogously, we can define the \emph{(strongly) } $(k,-)$\emph{-reachable} set $S_{\vec{G}, \sigma}^{k,-}(y)$ and the \emph{weakly } $(k,-)$\emph{-reachable} set $W_{\vec{G}, \sigma}^{k,-}(y)$  for any vertex $y \in V$; then define the $(k,-)$-\emph{coloring number} $\overrightarrow{\col}^-_{k}\left( \vec{G} \right) $  and the \emph{weak }$(k,-)$\emph{-coloring number}  $\overrightarrow{\wcol}^-_{k}\left( \vec{G} \right) $ of $\vec{G}$.

Define,  $ S_{\vec{G}, \sigma}^{k}\left( y\right)  :=  S_{\vec{G}, \sigma}^{k,+}\left( y\right) \bigcup   S_{\vec{G}, \sigma}^{k,-}\left( y\right)$  and$ W_{\vec{G}, \sigma}^{k}\left( y\right)  :=  W_{\vec{G}, \sigma}^{k,+}\left( y\right) \bigcup   W_{\vec{G}, \sigma}^{k,-}\left( y\right)$;  the vertices in $ S_{\vec{G}, \sigma}^{k}\left( y\right)$ (or $ W_{\vec{G}, \sigma}^{k}\left( y\right)$) are called the \emph{directed} \emph{(strongly)} $k$\emph{-reachable} (or directed weakly $k$-reachable) vertices from $y$.

Define the \emph{directed} $k$\emph{-coloring number} of $\vec{G}$, denoted by $\overrightarrow{\col}_{k}\left( \vec{G} \right) $, and the \emph{directed} \emph{weak }$k$\emph{-coloring number} 
of $\vec{G}$, denoted by $\overrightarrow{\wcol}_{k}\left( \vec{G}\right) $, by
\begin{eqnarray*}
\overrightarrow{\col}_{k}\left( \vec{G}\right) &=&\min_{\sigma\in \Pi }
\max_{y\in V\left( \vec{G} \right) }\left|  S_{\vec{G}, \sigma}^{k,+}(y) \bigcup S_{\vec{G}, \sigma}^{k,-}(y) \right|
\text{ and} \\
\overrightarrow{\wcol}_{k}\left(\vec{G}\right) &=&\min_{\sigma\in \Pi}
\max_{y\in V\left( \vec{G} \right) }\left|  W_{\vec{G}, \sigma}^{k,+}(y) \bigcup W_{\vec{G}, \sigma}^{k,-}(y) \right|;
\end{eqnarray*}

Let $G=(V,E)$ be a graph and $\mathcal{O}(G)$ be the set of all orientations of $G$.
Define the $(k,+)$-\emph{coloring number} of $\vec{G}$, denoted by $\overrightarrow{\col}^+_{k}\left( G\right) $,
and the \emph{weak} $(k,+)$\emph{-coloring number} of $G$, denoted by
$\overrightarrow{\wcol}^+_{k}\left( G\right) $, by
\begin{eqnarray*}
\overrightarrow{\col}^+_{k}\left( G\right) &=&\max_{\vec{G} \in \mathcal{O} \left( G\right)} \overrightarrow{\col}^+_{k}\left( \vec{G}\right)
\text{ and} \\\overrightarrow{\wcol}^+_{k}\left( G\right) &=&\max_{\vec{G} \in \mathcal{O} \left( G\right)} \overrightarrow{\wcol}^+_{k}\left( \vec{G}\right).
\end{eqnarray*}
Analogously, we can define the (\emph{strongly}) $(k,-)$\emph{-coloring number}  $\overrightarrow{\col}^-_{k}\left( {G} \right) $
and the \emph{weak }$(k,-)$\emph{-coloring number}  $\overrightarrow{\wcol}^-_{k}\left( {G} \right) $ of $G$.

Define the \emph{oriented} $k$-\emph{coloring number} of $G$, denoted by $\overrightarrow{\col}_{k}\left( G\right) $, and the \emph{oriented weak}
$k$\emph{-coloring number} of $G$, denoted by $\overrightarrow{\wcol}_{k}\left( G\right) $, by
\begin{eqnarray*}
\overrightarrow{\col}_{k}\left( G\right) &=&\max_{\vec{G} \in \mathcal{O} \left( G\right)} \overrightarrow{\col}_{k}\left( \vec{G}\right)
\text{ and} \\
\overrightarrow{\wcol}_{k}\left( G\right) &=&\max_{\vec{G} \in \mathcal{O} \left( G\right)} \overrightarrow{\wcol}_{k}\left( \vec{G}\right).
\end{eqnarray*}

Note that the ${\col}\left( G\right) =\overrightarrow{\col}_{1}\left(G\right) = \overrightarrow{\wcol}_{1}\left( G\right) $, $\overrightarrow{\col}_{k}\left(G\right) \leq {\col}_{k}\left( G\right) $, and $\overrightarrow{\wcol}_{k}\left(G\right) \leq {\wcol}_{k}\left( G\right) $.

When $k \geq \left|V\right| $ we refer to the \emph{complete oriented coloring number}, $\overrightarrow{\col}_{\infty }\left( G\right) $, and the \emph{complete oriented weak coloring number}, $\overrightarrow{\wcol}_{\infty }\left( G\right) $.

For the relationship between $\overrightarrow{\col}^+_{k}\left( G\right)$ and $\overrightarrow{\col}^-_{k}\left( G\right)$, we have the following observation:
\begin{observation} \label{CON-+}
For any graph $G$, any $k \ge 1$, we have $\overrightarrow{\col}^+_{k}\left( G\right) = \overrightarrow{\col}^-_{k}\left( G\right)$ and $\overrightarrow{\wcol}^+_{k}\left( G\right) = \overrightarrow{\wcol}^-_{k}\left( G\right)$.
\end{observation}

This observation follows from the following fact: for any orientation $\vec{G} \in \mathcal{O}(G)$, reversing all arcs of $\vec{G}$, produces the reverse orientation $\vec{G}^*$ of $\vec{G}$. This symmetry establishes Observation \ref{CON-+}.

The {\em maximum average degree} of a graph $G$, denoted by $\mad(G)$, is defined as
\[
\mad(G) = \max_{H \subseteq G} \frac {2|E(H)|} {|V(H)|},
\]
which places a bound on the average vertex degree in all subgraphs.
It is well-known (cf. \cite{SL-hak}, Theorem~$4$)  that a graph $G$ has an orientation $\vec{G}$
such that $\Delta^{+}\left(\vec{G} \right) \le k$ if and only if ${\rm mad}(G) \leq 2k$.
We obtain the following proposition.

\begin{proposition} \label{de-col-1}
If $G$ is a graph with ${\rm mad}(G) \leq 2k$, then $ \overrightarrow{\col}^+_{1}\left( G\right) \le k+1$.
\end{proposition}

For the relationship between $\overrightarrow{\col}_{k}\left( G\right)$ and $\overrightarrow{\col}^+_{k}({G} )$,
we pose the following questions.

\begin{problem} \label{CON--dd-+-}
For any graph $G$, are there linear functions $f, g$, such that for any $k \ge 1$ we have
$\overrightarrow{\col}_{k}\left( G\right) \le f(\overrightarrow{\col}^+_{k}\left( G\right))$ and $\overrightarrow{\wcol}_{k}\left( G\right) \le g( \overrightarrow{\wcol}^+_{k}\left( G\right))?$
In particular, do the following inequalities hold:
$\overrightarrow{\col}_{k}\left( G\right) \le 2\overrightarrow{\col}^+_{k}\left( G\right)$ and $\overrightarrow{\wcol}_{k}\left( G\right) \le 2 \overrightarrow{\wcol}^+_{k}\left( G\right)?$
\end{problem}

\begin{problem} \label{half2}
For any graph $G$, any $k \ge 1$, is there  an orientation $\vec{G}$ of $G$ such that
  $$\overrightarrow{\col}^+_{k}\left( \vec{G}\right)  = \left\lceil \frac{ \overrightarrow{\col}_{k}\left( G\right)}{2} \right\rceil, ~ ~  \overrightarrow{\wcol}^+_{k}\left( \vec{G}\right)  = \left\lceil \frac{ \overrightarrow{\wcol}_{k}\left( G\right)}{2} \right\rceil?$$
\end{problem}

The following lemma relates $\overrightarrow{\wcol}_{k}\left(G\right)$ and $ \overrightarrow{\col}_{k}\left( G\right)$,
which is analogous to Lemma \ref{kc-L1} (Lemma $3$ in \cite{KY}, Proposition $4.8$ in  \cite{Book-N-O}).

\begin{lemma} \label{kc-L1-D}
Let $G$ be a graph and  $k$ be a positive integer, then   $$\overrightarrow{\wcol}^+_{k}\left(G\right) \leq \left( \overrightarrow{\col}^+_{k}\left( G\right) \right) ^{k}, ~~\overrightarrow{\wcol}_{k}\left(G\right) \leq \left( \overrightarrow{\col}_{k}\left( G\right) \right) ^{k}.$$
\end{lemma}

\begin{proof}
Fix an orientation $\vec G$, an order $\sigma$, for any $k\geq 1$, we will use induction on $k$ to show that  $\left|  W_{\vec{G}, \sigma}^{k,+}(y) \right| \leq \left|S_{\vec{G}, \sigma}^{k,+} \left( y\right) \right|^{k}$ and $\left|  W_{\vec{G}, \sigma}^{k}(y) \right| \leq \left|  S_{\vec{G}, \sigma}^{k}(y) \right| ^{k}$ for all $y\in V(\vec{G})$.

The case $k=1$ is immediate. We now consider the induction step.

For each $x\in  W_{\vec{G}, \sigma}^{k,+}\left( y\right) $ let $\vec{P}_{yx}$ be the
shortest directed $y\rightsquigarrow x$ path such that every vertex $w$ of $\vec{P}_{yx}$ satisfies $x\leq
_{\sigma}w$. Let $z$ be the first vertex on $\vec{P}_{yx}$ such that $z<_{\sigma}y$ and let $
i$ be the distance from $y$ to $z$. Then $z \in  S_{\vec{G}, \sigma}^{i,+} \left( y\right)
- S_{\vec{G}, \sigma}^{i-1,+}\left( y\right) $ and $x\in  W_{\vec{G}, \sigma}^{k-i,+} \left( z\right) $. It follows
that
\[
\begin{split}
\left|W_{\vec{G}, \sigma}^{k,+}\left( y\right) \right|
& \leq \sum_{i=1}^{k}\left|S_{\vec{G}, \sigma}^{i,+} \left( y\right)
- S_{\vec{G}, \sigma}^{i-1,+}\left( y\right) \right| \left| W_{\vec{G}, \sigma}^{k-i,+} \left( z\right)\right| \\
& \leq \left|S_{\vec{G}, \sigma}^{k,+} \left( y\right) \right|\cdot\left| S_{\vec{G}, \sigma}^{k-1,+} \left( z\right)\right|^{k-1}   \mbox{ (using the induction hypothesis)}\\
 & \leq  \left|S_{\vec{G}, \sigma}^{k,+} \left( y\right) \right|^{k}.
\end{split}
\]

Minimize over $\sigma$ subject to $\max\limits_{y\in V(G)}\left|S_{\vec{G}, \sigma}^{k,+} (y) \right|$ being minimum,
we obtain $\overrightarrow{\wcol}_{k}^{+}(\vec{G}) \leq \left(\overrightarrow{\col_{k}^{+}}(\vec{G})\right)^{k}$.
Since $\vec{G}$ was arbitrary, we have $\overrightarrow{\wcol}_{k}^{+}(G) \leq \left(\overrightarrow{\col_{k}^{+}}(G)\right)^{k}$.


Recall $ W_{\vec{G}, \sigma}^{k}\left( y\right)  :=  W_{\vec{G}, \sigma}^{k,+}\left( y\right) \bigcup   W_{\vec{G}, \sigma}^{k,-}\left( y\right)$,  and $ S_{\vec{G}, \sigma}^{k}\left( y\right)  :=  S_{\vec{G}, \sigma}^{k,+}\left( y\right) \bigcup   S_{\vec{G}, \sigma}^{k,-}\left( y\right)$.
For each $x\in  W_{\vec{G}, \sigma}^{k}\left( y\right)$ let $\bar{P}_{yx}$ be the
shortest directed $y\rightsquigarrow x$ or  $x \rightsquigarrow y$ path such that every vertex $w$ of $\bar{P}_{yx}$ satisfies $x\leq_{\sigma}w$. Let $z$ be the first vertex on $\bar{P}_{yx}$ such that $z<_{\sigma}y$ and let $
i$ be the distance from $y$ to $z$. Then $z \in  S_{\vec{G}, \sigma}^{i} \left( y\right)
- S_{\vec{G}, \sigma}^{i-1}\left( y\right) $ and $x\in  W_{\vec{G}, \sigma}^{k-i} \left( z\right) $. It follows
that
\[
\begin{split}
\left|W_{\vec{G}, \sigma}^{k}\left( y\right)  \right| & \leq \sum_{i=1}^{k}\left|S_{\vec{G}, \sigma}^{i} \left( y\right) - S_{\vec{G}, \sigma}^{i-1}\left( y\right) \right| \left| W_{\vec{G}, \sigma}^{k-i} \left( z\right)\right| \\
& \leq \left|S_{\vec{G}, \sigma}^{k} \left( y\right) \right| \cdot\left| S_{\vec{G}, \sigma}^{k-1} \left( z\right)\right|^{k-1}   \mbox{ (using the induction hypothesis)}\\
 & \leq  \left|S_{\vec{G}, \sigma}^{k} \left( y\right) \right| ^{k}.
\end{split}
\]

Minimize over $\sigma$ such that $\max\limits_{y\in V}\left|S_{\vec{G}, \sigma}^{k} (y) \right|$ being minimum
and then maximize over $\vec G$ to obtain the second inequality.
\end{proof}

\section{Bounding q-transitive oriented chromatic number by oriented weak 2q-coloring number of graphs}

\label{S-q-dis-oc}
We bound the $q$-transitive oriented chromatic number in terms of the oriented weak $2q$-coloring number. The proof is analogous to that of Theorem 3.1 in  \cite{LY2026}.

\begin{theorem} \label{o-Tran-main}
For each graph $G$ and integer $q \ge 1$, there exists a function  $c_q(k_q)$, such that $\chi^q_o(G) \le 2^{k_q}-1$, 
for some  $k_q \le ((2q+1) \cdot \overrightarrow{\wc}_{2q}(G))^{\overrightarrow{\wc}_{2q}(G)}$.
\end{theorem}

\begin{proof}
Let $G$ be a graph and $\overrightarrow{\wc}_{2q}(G)=k$.
Let $\vec{G}$ be an arbitrary orientation of $G$ without any directed cycle of length $\le 2q$.
Let $\sigma$ be a linear ordering of $V(G)$ that witnesses $\overrightarrow{\wc}_{2q}(\vec{G}) \le k$.
Note that for any vertex $y \in W_{\vec{G}, \sigma}^{q}(x)$ and $y \neq x$, we cannot have both $y \in W_{\vec{G}, \sigma}^{q, +}(x)$ and $y \in W_{\vec{G}, \sigma}^{q, -}(x)$; otherwise, $\vec{G}$ will contain a directed cycle of length $\le 2q$.

First, we color the vertices of $G$ sequentially in the linear ordering $\sigma$, so that each
vertex $x$ is colored by a color $f(x)$ different from
the colors of vertices in $W_{\vec{G},\sigma}^{2q}(x)-\{x\}$. 
Observe the following: for any two vertices $ y_1, y_2$, if  $y_1 \in W_{\vec{G}, \sigma}^{q,+}(x)$, $y_2 \in W_{\vec{G}, \sigma}^{q,-}(x)$,  then the directed walk obtained by concatenating the paths witnessing paths of $y_1, y_2 \in W_{\sigma}^{q}(x)$ will contain a directed path that witnesses $y_1 \in W_{\vec{G}, \sigma}^{ 2q}(y_2)$ (if $y_1<_\sigma y_2$) or $y_2 \in W_{\vec{G}, \sigma}^{2q}(y_1)$ (if $y_2<_\sigma y_1$);
in either case, $f(y_2)\neq f(y_1)$.

In the above coloring process, when a vertex $x$ is colored by $f(x)$, suppose that
\[
\begin{split}
Y[x]  := ~ & \{y_i : 1 \le i \le t, ~ y_i  \in W_{\vec{G}, \sigma}^{q}(x), ~  y_1 <_\sigma y_2 <_\sigma \dots <_\sigma y_t = x \}.
\end{split}
\]

 For each $y_i\in Y[x]$, let the minimum length among all  the  {\em directed paths} that witness $y_i  \in W_{\sigma}^{q}(x) \cap Y[x]$ be $s_i$. If $y_i  \in W_{\vec{G},\sigma}^{q, +}(x)$, define $r_i := +s_i$; if $y_i  \in W_{\vec{G},\sigma}^{q, -}(x)$, define $r_i := -s_i$.  By the definition of $Y[x]$, $r_i$ is well defined for $y_i$, and $0 \le s_i \le q$.

For each color $\beta$ that has been used in the above coloring process for $f$, we split $\beta$ into $2q+1$ variants, denoted by $(\beta,0), (\beta, \pm{1}), \ldots, (\beta, \pm{q})$.
We then assign a label $\lambda(x)$ to $x$, where $\lambda(x):=  ((f(y_1),r_1), \ldots, (f(y_t),r_t))$.
Note that we always have  $(f(y_t),r_t) = (f(x),0)$ in $\lambda(x)$.
For counting purposes, suppose that, for all $x$,  we have $t \le t_q$; and $t_q \le \overrightarrow{\wc}_{2q}(\vec{G})$; indeed, if $\sigma$ also witnesses $\overrightarrow{\wc}_{q}(\vec{G})$, then $t_q \le \overrightarrow{\wc}_{q}(\vec{G})$.

For the induction step $q\ge2$, for orientation $\vec{G}$ of $G$,  we select vertices for coloring according to the linear ordering $\sigma$ of $V(G)$. Suppose we are about to color vertex $u$.
For each vertex to be colored $u$, we preassign $u$ a color $\lambda(u)$.

\begin{claim} \label{key-claim1}
For any two vertices $x$, $y$ in graph $\vec{G}$, if there is a directed path of length at most $2q$ between $x$ and $y$ in $\vec{G}$, then $\lambda(x) \neq \lambda(y)$.
\end{claim}

\begin{proof}
Suppose, for a contradiction, that there exist $x$, $y$ such that $\lambda(x_{0})=\lambda(x_{s})$ and there exists a directed path between $x$
and $y$ of length at most $2q$ in $\vec{G}$.
Suppose $\vec{P}=x_{0}(=x)x_{1}\cdots x_{s}(=y)$ is the shortest directed path between $x_{0}$ and $x_{s}$ in $\vec{G}$, and $s\leq 2q$, let $x_{l}$ be the $\sigma$-least vertex in $\vec{P}$.
If $x_{l}=x_{0}$ (or $x_{l}=x_{s}$),
then $x_{0}\in W^{2q}_{\vec{G},\sigma}(x_{s})$ (or $x_{s}\in W^{2q}_{\vec{G},\sigma}(x_{0})$), this implies that $f(x_{0}) \neq f(x_{s})$, this contradicts the assumption that $\lambda(x_{0})=\lambda(x_{s})$.

Otherwise, let $\vec{P}=\vec{P^{\prime}}+\vec{P^{\prime\prime}}=x_{0}\dots x_{l}+x_{l}x_{l+1}\dots x_{s}$ and $||\vec{P}||=||\vec{P^{\prime}}||+||\vec{P^{\prime\prime}}||\leq 2q$, note that min$\{||\vec{P^{\prime}}||, ||\vec{P^{\prime\prime}}||\}\leq q$.
Without loss of generality, assume that $||\vec{P^{\prime}}||\leq q$, then $\vec{P^{\prime\prime}}$ contains a directed path  that witnesses that $x_{l}\in W^{q,+}_{\vec{G},\sigma}(x_{0})$,
then $(f(x_{l}), r_{l})\in \lambda(x_{0})$, where $r_{l}=+s_{l}\leq ||\vec{P^{\prime}}||$.
Now we  claim that $(f(x_{l}), r_{l})\notin \lambda(x_{s})$.
On the one hand, $x_{l}\notin W^{s_{l}, +}_{\vec{G},\sigma}(x_{s})$, otherwise, there exists a directed path from $x_{s}$ to $x_{l}$ has length at most $|r_{l}|$ which witnesses that $x_{l}\in W^{s_{l}, +}_{\vec{G},\sigma}(x_{s})$, then the union of this directed path and $\vec{P^{\prime\prime}}$  contains a directed cycle of length $\leq 2q$, contradiction.
On the other hand, if there is a vertex $x_{l'}\neq x_{l}$, such that $(f(x_{l'}), r_{l'})=(f(x_{l}), r_{l})\in \lambda(x_{s})$,
this means that there exists a directed path that witnesses  $x_{l'}\in W^{|r_{l'}|, +}_{\vec{G}, \sigma}(x_{s})$ in $\vec{G}$ and the length of the directed path is $|r_{l'}|=|r_{l}|=+s_{l}\leq ||\vec{P^{\prime}}||$.
Denote this directed path as $\vec{P^{\prime\prime\prime}}=x_{s}\dots x_{l'}$, $||\vec{P^{\prime\prime\prime}}||=+s_{l}\leq ||\vec{P^{\prime}}||$.
Note that $V(\vec{P^{\prime\prime}})\cap V(\vec{P^{\prime\prime\prime}})=\emptyset$, so the union of $\vec{P^{\prime\prime}}$ and $\vec{P^{\prime\prime\prime}}$ contains a directed path that witnesses
$x_{l}\in W^{2q}_{\vec{G},\sigma}(x_{l'})$ (if $x_{l}<_{\sigma}x_{l'}$) or $x_{l'}\in W^{2q}_{\vec{G},\sigma}(x_{l})$ (if $x_{l'}<_{\sigma}x_{l}$),
in either case there are $f(x_{l'})\neq f(x_{l})$, this contradicts the assumption of $\lambda(x_{0})=\lambda(x_{s})$.
This proves the claim.
\end{proof}

For an orientation $\vec{G}$ of $G$, ${\vec{G}}^{q}$ is the $q$-th power of $\vec{G}$,
we let $U_{{\vec{G}}^{q}}$ be the undirected graph defined by $V(U_{\vec{G}^{q}})=V(\vec{G})$,
and $xy$ is an edge in $U_{{\vec{G}}^{q}}$ if and only if $x$ and $y$ are joined in $\vec{G}^{q}$ by a directed path of length at most $2$.
By Claim \ref{key-claim1}, $\lambda$ is a proper coloring of $U_{\vec{G}^{q}}$ and this assignment $\lambda$ uses at most $k_q = k_q( \overrightarrow{\wc}_{2q}(G)) :=  ((2q+1) \cdot \overrightarrow{\wc}_{2q}(G))^{t_q} \le ((2q+1) \cdot \overrightarrow{\wc}_{2q}(G))^{\overrightarrow{\wc}_{2q}(G)}$ colors.

As in the proof of Proposition \ref{pro-Pq-to-Oc}, let $\vec{H}_{k_{q}}$ be the oriented graph defined as follows: the vertices of $\vec{H}_{k_{q}}$ are all tuples of the form $(i; a_{1},\cdots, a_{i-1})$ with $1\leq i\leq k_{q}$ ($i$ is called the identity of the vertex) and $a_{j}\in \{0,1\}$ for every $j$, $1\leq j<i$.
Let $a=(i; a_{1},\cdots, a_{i-1})$ and $b=(j; b_{1},\cdots, b_{j-1})$ be two vertices in $\vec{H}_{k_{q}}$ with $i< j$.
There is an arc from $a$ to $b$ if $b_{i}=0$ and an arc from $b$ to $a$ otherwise.
Observe that the graph $\vec{H}_{k_{q}}$ has $2^{k_{q}}-1$ vertices and all arcs in $\vec{H}_{k_{q}}$ link vertices with distinct identities.

We denote the set of vertices $x$ with $\lambda(x)=i$ by $V_{i}$ and the (oriented) subgraph of $\vec{G}^{q}$ induced by $V_{i}\cup V_{j}$ by $\vec{G}^{q}_{i,j}$.
From the definition of $U_{{\vec{G}}^{q}}$ we get that ${\vec{G}}^{q}_{i, j}$ contains no directed 2-path in ${\vec{G}^{q}}$, i.e.,
if $x\in V_{i}$ and $x\rightarrow V_{j}$, then there is no vertex $y\in V_{j}$ such that $y\rightarrow x$.
We now define a map $\psi: V({\vec{G}}^{q})$ to $V(\vec{H}_{k_{q}})$, and for every $x\in V_{i}, 1\leq i\leq s$, let $\psi(x)=(i; a_{1},\cdots, a_{i-1})$ where for every $j$, $1\leq j< i$, $a_{j}=0$ if there is an arc from $V_{j}$ to $x$ in ${\vec{G}}^{q}$ and $a_{j}=1$ otherwise. This implies that $\psi$ is an oriented $\vec{H}_{k_{q}}$-coloring of ${\vec{G}}^{q}$.
Thus $\psi$ is clearly a $q$-transitive $\vec{H}_{k_{q}}$-coloring of $\vec{G}$. This finishes the proof.
\end{proof}

By Proposition \ref{pro-oc-Pq} and Theorem \ref{o-Tran-main}, we have the following corollary.

\begin{corollary} \label{cor-main-1}
Suppose a graph $G$ has $d_q: = \overrightarrow{\wc}_{2q}(G)$,  then  for any orientation  $\vec{G}$ of $G$ without directed cycles of length $2q$, there exists a function $c_q(d_q)$, such that $\vec{\chi}_{\vec{P}_{2q}} (\vec{G}) \le c_q(d_q)$.
\end{corollary}

Setting $q=1$ in Theorem \ref{o-Tran-main}, we obtain $\chi_o(G) \le c_1(k_1)$ and
$k_1 \le (3 \overrightarrow{\wc}_{2}(G))^{\overrightarrow{\wc}_{2}(G)}$ immediately. However, the case $q=1$ admits a simpler treatment than the general case; by adapting the proof of Theorem \ref{o-Tran-main} to $q=1$, we obtain the following theorem.

\begin{theorem} \label{o-co-2}
There exists a function $\bar{c}_1(k)$, such that for any graph $G$, $\chi_o(G) \le \bar{c}_1(\overrightarrow{\wc}_2(G))$;
indeed, $\chi_o(G) \le \overrightarrow{\wc}_2(G)\cdot2^{\overrightarrow{\wc}_2(G)-1}$.
\end{theorem}

For comparison with existing results, consider the following theorem of Aravind and Subramanian in \cite{AS}.

\begin{theorem} [\cite{AS}]
If $G$ is an undirected graph with maximum degree $\Delta$ and degeneracy $d$, then
$\chi_o(G) \le 16 \cdot \Delta \cdot d \cdot 2^d$.
\end{theorem}

By Proposition \ref{pro-oc-Pq} and Theorem \ref{o-co-2}, we have the following corollary.

\begin{corollary} \label{cor-co-2}
Suppose a graph $G$ has $p: = \overrightarrow{\wc}_2(G)$, then  for any orientation  $\vec{G}$ of $G$, we have $\vec{\chi}_{\vec{P}_{2}} (\vec{G}) \le p\cdot2^{p-1}$.
\end{corollary}

\section {Bounding generalized coloring numbers by using oriented generalized coloring numbers of graphs}

\label{sec-re-to-co}

In this section, we first consider a graph $H$ with a subdivision model in $G$ in which each replacement path has length at most $l$. We show that a bound on this parameter for $G$ at radius $l$ gives a bound on the degeneracy of $H$. We then apply standard inequalities relating admissibility and shallow topological density. This gives an  upper bound on the generalized coloring number of $G$ in terms of its oriented generalized coloring number.

\begin{definition}
Let $t\ge 0$ be an integer.
A graph $H$ is a \emph{depth-$t$ shallow topological minor} of $G$ if a subdivision of $H$ occurs as a subgraph of $G$,
with each edge of $H$ replaced by a path of length at most $2t+1$.
We use the notation $H \preccurlyeq^{top}_{t} G$ to indicate that $H$ is a shallow topological minor of $G$ at depth $t$.

The topological greatest reduced average density (top-grad) of rank $r$ of a graph $G$ is defined by:
$$\widetilde{\nabla}_{t}(G)=\max\left\{\frac{|E(H)|}{|V(H)|}: H \preccurlyeq^{top}_{t} G \right\}.$$
\end{definition}

In a subdivision model, the vertices corresponding to $V(H)$ are called \emph{branch vertices}.
The replacement paths have pairwise disjoint interiors, and no branch vertex is an internal vertex of a replacement path.

The following lemma provides the key connection between directed generalized coloring numbers and top-grad.
Recall that $\deg(H)$ is the least integer $d$ such that every nonempty subgraph of $H$ contains a vertex of degree at most $d$.

\begin{lemma}
\label{lem:short-subdivision}
Let $l,~p\ge 1$  and suppose that $\overrightarrow{\col}_{l}(G)\le p$.
If  $H$ has a subdivision model in $G$ in which every replacement path has length at most $l$, then $\deg(H)\le p-1$.\\
Consequently,
\[
\frac{|E(H)|}{|V(H)|}\leq p-1.
\]
\end{lemma}
\begin{proof}
Fix such a subdivision model and identify each vertex of $H$ with its corresponding branch vertex in $G$.
For each edge $xy\in E(H)$, let $P_{xy}$ be its replacement path.

Orient each replacement path consistently in either direction, so that it becomes a directed path.
Since distinct replacement paths are edge-disjoint, these choices are compatible.
Extend them arbitrarily to an orientation $\vec G$ of all of $G$.

Since $\overrightarrow{\col}_{\ell}(\vec{G})\le \overrightarrow{\col}_{\ell}(G)\le p$, there is an order $\sigma\in\Pi(G)$ such that
\[
   |S_{\vec{G},\sigma}^{\ell}(y)|\le p
   \qquad\text{for every }y\in V(G).
\]

Let $H'$ be an arbitrary nonempty subgraph of $H$. Choose $v\in V(H')$ to be the branch vertex occurring last in $\sigma$ among the vertices of $H'$. Every $u\in N_{H'}(v)$ therefore satisfies $u<_\sigma v$.

For each $u\in N_{H'}(v)$, let $x_u$ be the first vertex on $P_{vu}$ from $v$ to $u$ that precedes $v$ in $\sigma$.
 Such a vertex exists because $u<_\sigma v$.
By construction, the subpath from $v$ to $x_u$ has length at most $l$, and all its internal vertices $w$ satisfy $v<_\sigma w$.
Moreover, this subpath is directed either from $v$ to $x_u$ or from $x_u$ to $v$. Hence $x_{u}$ is directed $l$-reachable from $v$.

Distinct neighbors $u$ of $v$ in $H'$ give rise to distinct vertices $x_u$.
Indeed, the replacement paths associated with distinct edges incident with $v$ intersect only at $v$, and none of the vertices $x_u$ equals $v$.
Since $v$ is directed $l$-reachable from itself, we obtain
\[
d_{H'}(v)
\le
\left|S_{\vec{G},\sigma}^{\ell}(y)\right|\le -1
\le p-1.
\]
Thus every nonempty subgraph of $H$ contains a vertex of degree at most $p-1$, proving
$\deg(H)\le p-1$.

Finally, successively deleting a vertex of degree at most $p-1$ counts each edge exactly once and gives
\[
\frac{|E(H)|}{|V(H)|}\leq p-1.
\]
\end{proof}

\begin{corollary}
\label{cor:density}
Let $t\ge 0$ be an integer. If $\overrightarrow{\col}_{2t+1}(G)\le p$, then $\widetilde{\nabla}_{t}(G)\le p-1$.
\end{corollary}

\begin{proof}
Let $H$ be a depth-$t$ shallow topological minor of $G$.
By definition, $H$ has a subdivision model in $G$ in which each edge of $H$ is represented by a path of length at most $2t+1$.
By Lemma~\ref{lem:short-subdivision},
\[
  \frac{|E(H)|}{|V(H)|} \le p-1.
\]
 This completes the proof.
\end{proof}

We use the following standard estimates, in a form sufficient
for our purposes; see \cite{S.S2026}.

\begin{proposition}
\label{prop:standard}
For each graph $G$ and integer $r\ge 1$,
$$\col_r(G) \le \wcol_r(G) \le (6r)^r \cdot \widetilde{\nabla}_{r-1}(G)^{3r}.$$
\end{proposition}


We can now combine Corollary~\ref{cor:density} with the standard admissibility estimates.

\begin{theorem}
\label{M-g-d}
Let $G$ be a  graph, $k\ge 1$ be an integer, and let $p\ge 1$.
If $\overrightarrow{\col}_{2k-1}(G)\le p$, then $\col_k(G)\le (6k)^{k}(p-1)^{3k}.\label{eq:main-bound}$
\end{theorem}

\begin{proof}
Applying Corollary~\ref{cor:density} with $t=k-1$ gives $\widetilde{\nabla}_{k-1}(G)\le p-1$.
By Proposition~\ref{prop:standard}, we conclude that
$$\col_{k}(G) \leq(\adm_{k}(G))^{k} \leq (6k)^{k}(p-1)^{3k}.$$
\end{proof}

By \cite[Corollary~4.10]{KRSW}, a digraph class has directed bounded expansion if and only if its directed weak $r$-coloring numbers are uniformly bounded for each fixed $r\ge 1$. We apply this characterization to the class of all orientations of graphs in $\mathcal C$.
For this class, the condition is equivalent to a uniform bound on $\overrightarrow{\mathrm{wcol}}_r(G)$ over all $G\in\mathcal C$ for each fixed $r$.
Since $\overrightarrow{\mathrm{wcol}}_r(G)\le \mathrm{wcol}_r(G)$, bounded expansion of $\mathcal C$ gives such a bound.
Conversely, $\overrightarrow{\mathrm{col}}_r(G)\le \overrightarrow{\mathrm{wcol}}_r(G)$, so Theorem~\ref{M-g-d} gives a uniform bound on $\mathrm{col}_k(G)$ by taking $r=2k-1$.
This proves the following corollary.
\begin{corollary}\label{undirected-directed-bounded-expension}
An undirected graph class $\mathcal{C}$ has bounded expansion if and only if the digraph class $\mathcal{D}=\{\vec{G}: \vec{G} \in \mathcal{O}(G)$ for each $G\in\mathcal{C}\}$ has directed bounded expansion.
\end{corollary}

\section{Conclusions}
In this paper, we introduced the $q$-transitive oriented chromatic number and bounded it in terms of the oriented weak $2q$-coloring number.
We also compared $q$-transitive oriented colorings with $\vec{P}_{2q}$-preserving oriented colorings, and showed that a graph class has bounded generalized
coloring numbers  if and only if it has bounded oriented generalized coloring numbers.
The following theorem summarizes the resulting characterizations of bounded expansion.

\begin{theorem} \label{T-list}
	Let $\Phi$ be a class of graphs,  the following statements  are equivalent.
	\begin{enumerate}
		\item $\Phi$ has bounded expansion;
		
		\item For every integer $k$, there is a constant $a_k$ such that $\vec{\chi}_{\vec{P}_{2k}} ({G}) \le a_k$, for all $G \in \Phi$ and every orientation $\vec G$ of $G$ with no directed cycle of length at most $2k$;
		
		\item For every integer $k$, there is a constant $b_k$ such that $\chi^k_o(G) \le b_k$, for all $G \in \Phi$;

		\item For every integer $k$, there is a target oriented graph $\vec{H}_k$ with no directed cycle of length at most $k$, such that, for any acyclic orientation $\vec{G}$ of $G \in \Phi$, there exists an oriented homomorphism from $\vec{G}$ to  $\vec{H}_k$;
		
		\item For every integer $k$, there is a constant $c_k$ such that the weak $k$-coloring number ${\wcol}_{k}\left({G} \right)  \le c_k$, for all $G \in \Phi$.
		
		\item For every integer $k$, there is a constant $d_k$ such that the oriented weak $k$-coloring number $\overrightarrow{\wcol}_{k}\left({G} \right)  \le d_k$, for all $G \in \Phi$.
	\end{enumerate}
\end{theorem}

\begin{proof}	
	By Theorem 11.4 in \cite{Book-N-O}, we have, in Theorem \ref{T-list}, Statement 1 and Statement 4 are equivalent.
	By Proposition \ref{pro-oc-Pq} and Proposition \ref{pro-Pq-to-Oc} of this paper, we have Statement 2 and Statement 3 are equivalent.
	The equivalence between Statement 1 and Statement 5 was proved  by Zhu in \cite{Zhu-BGCN} and  Yang in \cite{Yang-GTF}.
	The equivalence between Statement 5 and Statement 6 follows from the definitions and Theorem \ref{M-g-d} of this paper.
	By Theorem \ref{o-Tran-main} of this paper, Statement 6 $\Rightarrow$ Statement 3. Finally, by definition, we have Statement 3 $\Rightarrow$ Statement 4.
	This finishes the proof of Theorem \ref{T-list}. 
\end{proof}

\bibliographystyle{plain}
\bibliography{References}

\end{document}